\documentclass[12pt]{article}
\usepackage{amsfonts}
\usepackage{amssymb}
\title{Vertex operator algebras
associated to certain admissible modules for affine Lie algebras of
type $A$}
\author{Ozren Per\v{s}e}
\date{}
\begin{document}
\def \Z{\Bbb Z}
\def \C{\Bbb C}
\def \R{\Bbb R}
\def \Q{\Bbb Q}
\def \N{\Bbb N}
\def \tr{{\rm tr}}
\def \span{{\rm span}}
\def \Res{{\rm Res}}
\def \End{{\rm End}}
\def \E{{\rm End}}
\def \Ind {{\rm Ind}}
\def \Irr {{\rm Irr}}
\def \Aut{{\rm Aut}}
\def \Hom{{\rm Hom}}
\def \mod{{\rm mod}}
\def \ann{{\rm Ann}}
\def \<{\langle}
\def \>{\rangle}
\def \t{\tau }
\def \a{\alpha }
\def \e{\epsilon }
\def \l{\lambda }
\def \L{\Lambda }
\def \g{\gamma}
\def \b{\beta }
\def \om{\omega }
\def \o{\omega }
\def \c{\chi}
\def \ch{\chi}
\def \cg{\chi_g}
\def \ag{\alpha_g}
\def \ah{\alpha_h}
\def \ph{\psi_h}
\def \be{\begin{equation}\label}
\def \ee{\end{equation}}
\def \bl{\begin{lem}\label}
\def \el{\end{lem}}
\def \bt{\begin{thm}\label}
\def \et{\end{thm}}
\def \bp{\begin{prop}\label}
\def \ep{\end{prop}}
\def \br{\begin{rem}\label}
\def \er{\end{rem}}
\def \bc{\begin{coro}\label}
\def \ec{\end{coro}}
\def \bd{\begin{de}\label}
\def \ed{\end{de}}
\def \pf{{\bf Proof. }}
\def \voa{{vertex operator algebra}}

\newtheorem{thm}{Theorem}[section]
\newtheorem{prop}[thm]{Proposition}
\newtheorem{coro}[thm]{Corollary}
\newtheorem{conj}[thm]{Conjecture}
\newtheorem{lem}[thm]{Lemma}
\newtheorem{rem}[thm]{Remark}
\newtheorem{de}[thm]{Definition}
\newtheorem{hy}[thm]{Hypothesis}
\makeatletter \@addtoreset{equation}{section}
\def\theequation{\thesection.\arabic{equation}}
\makeatother \makeatletter

\newcommand{\binom}[2]{{{#1}\choose {#2}}}
    \newcommand{\nno}{\nonumber}
    \newcommand{\lbar}{\bigg\vert}
    \newcommand{\p}{\partial}
    \newcommand{\dps}{\displaystyle}
    \newcommand{\bra}{\langle}
    \newcommand{\ket}{\rangle}
 \newcommand{\res}{\mbox{\rm Res}}
\renewcommand{\hom}{\mbox{\rm Hom}}
  \newcommand{\epf}{\hspace{2em}$\Box$}
 \newcommand{\epfv}{\hspace{1em}$\Box$\vspace{1em}}
\newcommand{\nord}{\mbox{\scriptsize ${\circ\atop\circ}$}}
\newcommand{\wt}{\mbox{\rm wt}\ }

\maketitle
\begin{abstract}
Let $L(-\frac{1}{2}(l+1),0)$ be the simple vertex operator algebra
associated to an affine Lie algebra of type $A_{l}^{(1)}$ with the
lowest admissible half-integer level $-\frac{1}{2}(l+1)$, for
even~$l$. We study the category of weak modules for that vertex
operator algebra which are in category $\cal{O}$ as modules for the
associated affine Lie algebra. We classify irreducible objects in
that category and prove semisimplicity of that category.
\end{abstract}

\footnotetext[2]{
Partially supported by the
Ministry of Science, Education and Sports
of the Republic of Croatia, grant 0037125.}


\section{Introduction}

Let ${\frak g}$ be a simple finite-dimensional Lie algebra and
$\hat{\frak g}$ the associated affine Lie algebra. For any complex
number $k\ne -h^{\vee}$, denote by $L(k,0)$ the simple vertex
operator algebra associated to $\hat{\frak g}$ with level $k$. The
representation theory of $L(k,0)$ heavily depends on the choice of
level $k$. If $k$ is a positive integer, $L(k,0)$ is a rational
vertex operator algebra (cf. \cite{FZ}, \cite{Z}), i.e. the category
of $\Z _{+}$-graded weak $L(k,0)$-modules is semisimple. Irreducible
objects in that category are integrable highest weight $\hat{\frak
g}$-modules of level $k$ (\cite{FZ}, \cite{L}). The corresponding
associative algebra $A(L(k,0))$, defined in \cite{Z}, is
finite-dimensional (cf. \cite{KWn}). In some cases such as $k \notin
\Q$ or $k < -h^{\vee}$ (studied in \cite{KL1} and \cite{KL2}),
categories of $L(k,0)$-modules have significantly different
structure then categories of $L(k,0)$-modules for a positive integer
$k$. However, there are examples of rational levels $k$ such that
the category of weak $L(k,0)$-modules which are in the category
$\cal{O}$ as $\hat{\frak g}$-modules, has similar structure as the
category of $\Z _{+}$-graded weak $L(k,0)$-modules for positive
integer levels $k$. These are so called admissible levels, defined
in \cite{KW1} and \cite{KW2}.

The case of vertex operator algebras associated to affine Lie
algebras of type $C_{l}^{(1)}$ with admissible half-integer levels
has been studied in \cite{A1} and \cite{A2}. Vertex operator
algebras associated to affine Lie algebras of type $A_{1}^{(1)}$
with arbitrary admissible level have been studied in \cite{AM} and
\cite{DLM}. In these cases vertex operator algebra $L(k,0)$ has
finitely many irreducible weak modules from the category $\cal{O}$
and every weak $L(k,0)$-module from the category $\cal{O}$ is
completely reducible. One can say that these vertex operator
algebras are rational in the category $\cal{O}$. In \cite{AM},
authors gave a conjecture that vertex operator algebras $L(k,0)$,
for all admissible levels $k$, are rational in the category
$\cal{O}$. In the case of vertex operator algebras associated to
affine Lie algebras of type $B_{l}^{(1)}$ with admissible
half-integer levels, certain parts of this conjecture were verified
in \cite{P1}. Admissible modules for affine Lie algebras were also
recently studied in \cite{A3}, \cite{FM}, \cite{GPW}, \cite{P2},
\cite{W}. Vertex operator algebras associated to certain affine Lie
algebras with non-admissible negative integer levels have been
studied in \cite{AP}.

When $k$ is an admissible level, vertex operator algebra $L(k,0)$ is
a quotient of the generalized Verma module by the maximal ideal
generated by one singular vector. The formula for this singular
vector is very complicated for general admissible level $k$ (cf.
\cite{MFF}). But for some special cases of affine Lie algebras and
half-integer admissible levels $k$, this singular vector has
conformal weight $2$, and the formula for this vector is relatively
simple (cf. \cite{A1}, \cite{A2} and \cite{P1}). In this paper we
study one similar special case, for which we verify the conjecture
from \cite{AM}.

We consider the case of an affine Lie algebra of type $A_{l}^{(1)}$
and the corresponding vertex operator algebra
$L(-\frac{1}{2}(l+1),0)$, for even $l$. We show that
$-\frac{1}{2}(l+1)$ is an admissible level for this affine Lie
algebra. The results on admissible modules from \cite{KW1} imply
that $L(-\frac{1}{2}(l+1),0)$ is a quotient of the generalized Verma
module by the maximal ideal generated by a singular vector of
conformal weight $2$. Using results from \cite{FZ}, \cite{Z}, we can
identify the corresponding associative algebra
$A(L(-\frac{1}{2}(l+1),0))$ with a certain quotient of $U({\frak
g})$. Algebra $A(L(-\frac{1}{2}(l+1),0))$ is infinite-dimensional in
this case. Using methods from  \cite{A2}, \cite{AM} and \cite{MP},
we get that irreducible $A(L(-\frac{1}{2}(l+1),0))$-modules from the
category $\mathcal{O}$ are in one-to-one correspondence with zeros
of the certain set of polynomials ${\mathcal P}_{0}$. By determining
a basis for the vector space ${\mathcal P}_{0}$, we obtain the
classification of irreducible $A(L(-\frac{1}{2}(l+1),0))$-modules
from the category $\mathcal{O}$. Using results from \cite{Z}, we
obtain the classification of irreducible weak
$L(-\frac{1}{2}(l+1),0)$-modules from the category $\mathcal{O}$.
Using this classification and results from \cite{KW2}, we show that
every weak $L(-\frac{1}{2}(l+1),0)$-module from the category
$\mathcal{O}$ is completely reducible.

In the case when $l$ is odd, the lowest half-integer admissible
level for affine Lie algebra of type $A_{l}^{(1)}$ is
$-\frac{1}{2}l$, and the maximal submodule of the generalized Verma
module is generated by a singular vector of conformal weight $4$. It
is more complicated to determine the formula for singular vector in
that case, and to use the method for classification from \cite{A2},
\cite{AM} and \cite{MP}.

In this paper ${\Z}_{+}$ denotes the set of nonnegative integers.

This paper (together with \cite{P1}) is a part of author's Ph.D.
thesis. The author would like to thank Professor Dra\v{z}en Adamovi\'{c} and
Professor Mirko Primc for their long time support and for numerous
helpful advice on this and related work.

\section{Vertex operator algebras associated to affine Lie algebras}

In this section we review certain results on vertex operator
algebras and corresponding modules. Specially, we recall some
results on vertex operator algebras associated to affine Lie
algebras.

\subsection{Vertex operator algebras and modules}

Let $(V, Y, {\bf 1}, \omega)$ be a vertex operator algebra (cf.
\cite{B}, \cite{FHL} and \cite{FLM}). An {\it ideal} in a vertex
operator algebra $V$ is a subspace $I$ of $V$ satisfying $ Y(a,z)I
\subseteq I[[z,z^{-1}]] $ for any $a \in V$. Given an ideal $I$ in
$V$, such that ${\bf 1} \notin I$, $\omega \notin I$, the quotient
$V/I$ admits a natural vertex operator algebra structure.

Let $(M, Y_{M})$ be a weak module for a vertex operator algebra $V$
(cf. \cite{L}). A {\it ${\Z}_{+}$-graded weak} $V$-module
(\cite{FZ}) is a weak $V$-module $M$ together with a
${\Z}_{+}$-gradation $M=\oplus_{n=0}^{\infty}M(n)$ such that
\begin{eqnarray*}
& &a_{m}M(n)\subseteq M(n+r-m-1)\;\;\;\mbox{for }a \in
V_{(r)},m,n,r\in {\Z},
\end{eqnarray*}
where $M(n)=0$ for $n < 0$ by definition.

A weak $V$-module $M$ is called a {\it $V$-module} if $L(0)$
acts semisimply on $M$ with the decomposition into
$L(0)$-eigenspaces $M=\oplus_{\alpha\in {\C}}M_{(\alpha)}$ such
that for any $\alpha \in {\C}$, $\dim M_{(\alpha)}<\infty$
and $M_{(\alpha+n)}=0$ for $n\in {\Z}$ sufficiently small.

\subsection{Zhu's $A(V)$ theory}

Let $V$ be a vertex operator algebra. Following \cite{Z}, we define
bilinear maps $* :  V \times V \to V$ and $\circ  :  V \times V \to
V$ as follows. For any homogeneous $a \in V$ and for any $b \in V$,
let
\begin{eqnarray}
a \circ b=\Res_{z}\frac{(1+z)^{{\rm wt} a}}{z^{2}}Y(a,z)b, \nonumber \\
a * b=\Res_{z}\frac{(1+z)^{{\rm wt} a}}{z}Y(a,z)b \nonumber
\end{eqnarray}
and extend to $V \times V \to V$ by linearity. Denote by
$O(V)$ the linear span of elements of the form $a \circ b$, and by
$A(V)$ the quotient space $V/O(V)$. For $a \in V$, denote by
$[a]$ the image of $a$ under the projection of $V$ onto
$A(V)$. The multiplication $*$ induces the multiplication
on $A(V)$ and $A(V)$ has a structure of an associative algebra.

\begin{prop}[\cite{FZ}, Proposition 1.4.2] Let $I$ be an ideal of $V$. Assume
${\bf 1} \notin I$, $\omega \notin I$. Then the associative
algebra $A(V/I)$ is isomorphic to $A(V)/A(I)$, where
$A(I)$ is the image of $I$ in $A(V)$.
\end{prop}

For any homogeneous $a \in V$ we define $o(a)=a_{{\rm wt} a-1}$
and extend this map linearly to $V$.

\begin{prop} [\cite{Z}, Theorem 2.1.2, Theorem 2.2.1]
\item[(a)] Let $M=\oplus_{n=0}^{\infty}M(n)$ be a ${\Z}_{+}$-graded
weak $V$-module. Then $M(0)$ is an $A(V)$-module defined as follows:
$$[a].v=o(a)v,$$
for any $a \in V$ and $v \in M(0)$.
\item[(b)] Let $U$ be an
$A(V)$-module. Then there exists a ${\Z}_{+}$-graded weak $V$-module
$M$ such that the $A(V)$-modules $M(0)$ and $U$ are isomorphic.
\end{prop}

\begin{prop}[\cite{Z}, Theorem 2.2.2] \label{t.1.3.5}
The equivalence classes of the irreducible $A(V)$-modules and the
equivalence classes of the irreducible ${\Z}_{+}$-graded weak
$V$-modules are in one-to-one correspondence.
\end{prop}

\subsection{Modules for affine Lie algebras}

Let ${\frak g}$ be a simple Lie algebra over ${\C}$ with a
triangular decomposition \linebreak ${\frak g}={\frak n}_{-} \oplus
{\frak h} \oplus {\frak n}_{+}$. Let $\Delta$ be the root system of
$({\frak g}, {\frak h})$, $\Delta_{+}\subset \Delta$
the set of positive roots, $\theta$ the highest root and
$(\cdot, \cdot): {\frak g}\times {\frak g}\to {\C}$ the Killing form,
normalized by the condition $(\theta, \theta)=2$.

The affine Lie algebra $\hat{\frak g}$ associated to ${\frak g}$ is
the vector space ${\frak g}\otimes {\C}[t, t^{-1}] \oplus {\C}c $
equipped with the usual bracket operation and the canonical central
element~$c$ (cf. \cite{K1}). Let $h^{\vee}$ be the dual Coxeter
number of $\hat{\frak g}$. Let $\hat{\frak g}=\hat{\frak n}_{-}
\oplus \hat{\frak h} \oplus \hat{\frak n}_{+}$ be the corresponding
triangular decomposition of $\hat{\frak g}$. Denote by
$\hat{\Delta}$ the set of roots of $\hat{\frak g}$, by
$\hat{\Delta}^{\mbox{\scriptsize{re}}}$ (resp.
$\hat{\Delta}^{\mbox{\scriptsize{re}}}_{+}$) the set of real (resp.
positive real) roots of $\hat{\frak g}$ and by $\alpha ^{\vee}$
denote the coroot of a real root $\alpha \in
\hat{\Delta}^{\mbox{\scriptsize{re}}}$.

For every weight $\lambda\in \hat{{\frak h}}^{*}$, denote by
$M( \lambda)$ the Verma module for $\hat{\frak g}$ with highest
weight $\lambda$, and by $L( \lambda)$ the irreducible
$\hat{\frak g}$-module with highest
weight~$\lambda$.

Let $U$ be a ${\frak g}$-module, and let $k\in {\C}$. Let
$\hat{\frak g}_{+}={\frak g}\otimes t{\C}[t]$ act trivially on $U$
and $c$ as scalar $k$. Considering $U$ as a ${\frak g}\oplus {\C}c
\oplus \hat{\frak g}_{+}$-module, we have the induced $\hat{\frak
g}$-module (so called {\it generalized Verma module})
$$N(k,U)=U(\hat{\frak g})\otimes_{U({\frak g}\oplus {\C}c
\oplus \hat{\frak g}_{+})} U.$$

For a fixed $\mu \in {\frak h}^{*}$, denote by $V(\mu)$ the
irreducible highest weight ${\frak g}$-module with highest weight
$\mu$. We shall use the notation $N(k, \mu)$ to denote the
$\hat{\frak g}$-module $N(k,V(\mu))$. Denote by $J(k, \mu)$ the
maximal proper submodule of $N(k, \mu)$ and  $L(k, \mu)=N(k, \mu)
/J(k, \mu)$.

\subsection{Admissible modules for affine Lie algebras}

Let $\hat{\Delta}^{\vee \mbox{\scriptsize{re}}}$\ (resp.
$\hat{\Delta}^{\vee \mbox{\scriptsize{re}}}_{+}$) $\subset
\hat{\frak h}$ be the set of real (resp. positive real) coroots of $
\hat{\frak g}$. Fix $\lambda \in \hat{\frak h}^*$. Let
$\hat{\Delta}^{\vee \mbox{\scriptsize{re}}}_{\lambda}= \{\alpha\in
\hat{\Delta}^{\vee \mbox{\scriptsize{re}}} \ |\
\langle\lambda,\alpha\rangle\in{\Z} \}$, $\hat{\Delta}^{\vee
\mbox{\scriptsize{re}}}_{\lambda +}=\hat{\Delta}^{\vee
\mbox{\scriptsize{re}}}_{\lambda} \cap \hat{\Delta}^{\vee
\mbox{\scriptsize{re}}}_{+}$, $\hat{\Pi}^{\vee}$ the set of simple
coroots in $\hat{\Delta}^{\vee \mbox{\scriptsize{re}}}$ and
$\hat{\Pi}^{\vee}_{\lambda }= \{\alpha \in \hat{\Delta}^{\vee
\mbox{\scriptsize{re}}}_{\lambda +}\ \vert \ \alpha$ not equal to a
sum of several coroots from $\hat{\Delta}^{\vee
\mbox{\scriptsize{re}}}_{\lambda +} \}$. Define $\rho$ in the usual
way, and denote by $w. \lambda$ the "shifted" action of an element
$w$ of the Weyl group of $ \hat{\frak g}$.

Recall that a weight  $\lambda \in \hat{\frak h} ^* $ is called {\it
admissible} (cf. \cite{KW1}, \cite{KW2} and \cite{W}) if the
following properties are satisfied:
\begin{eqnarray*}
& &\langle\lambda + \rho,\alpha\rangle \notin -{\Z}_+ \mbox{ for all }
\alpha \in \hat{\Delta}^{\vee \mbox{\scriptsize{re}}}_{+}, \\
& &{\Q} \hat{\Delta}^{\vee \mbox{\scriptsize{re}}}_{\lambda}={\Q} \hat{\Pi}^{\vee} .
\end{eqnarray*}
The irreducible $\hat{\frak g}$-module $L(\lambda)$
is called admissible if the weight $\lambda \in \hat{\frak h} ^* $
is admissible.

We shall use the following results from \cite{KW1} and \cite{KW2}:
\begin{prop}[\cite{KW1}, Corollary 2.1] \label{t.KW1}
Let $\lambda $ be an admissible weight. Then
$$L(\lambda)=\frac {M(\lambda)}{\sum_{\alpha \in \hat{\Pi}^{\vee}_{\lambda }}U (
\hat{\frak g}) v^{\alpha}}\ ,$$
where $v^{\alpha}\in M(\lambda)$ is a singular vector of weight $r_{\alpha}.
\lambda$, the highest weight vector of $M(r_{\alpha}.\lambda)=\ U(\hat{\frak g})
v^{\alpha}\subset M(\lambda)$.
\end{prop}

\begin{prop}[\cite{KW2}, Theorem 4.1] \label{t.KW2}
Let $M$ be a $\hat{\frak g}$-module from the category $\mathcal{O}$ such
that for any irreducible subquotient $L(\nu )$ the weight $\nu $ is
admissible. Then $M$ is completely reducible.
\end{prop}

\subsection{Vertex operator algebras $N(k,0)$ and $L(k,0)$, for
 $k \neq - h^{\vee}$}

Since $V(0)$ is the one-dimensional
trivial ${\frak g}$-module, it can be identified with~${\C}$.
Denote by ${\bf 1}=1 \otimes 1 \in N(k,0)$.
We note that $N(k,0)$ is spanned by the elements of the form
$x_{1}(-n_{1}-1)\cdots x_{m}(-n_{m}-1){\bf 1}$, where  $x_{1}, \dots, x_{m}\in
{\frak g}$ and $n_{1}, \dots, n_{m}\in {\Z}_{+}$, with $x(n)$ denoting
the representation image of
$x\otimes t^{n}$ for $x\in{\frak g}$ and $n\in {\Z}$.
Vertex operator map
$Y(\cdot, z): N(k,0) \to  (\mbox{\rm End}\;N(k,0))[[z, z^{-1}]]$
is uniquely determined by defining $Y({\bf 1}, z)$ to be the
identity operator on $N(k,0)$ and
$$Y(x(-1){\bf 1}, z)=\sum_{n\in {\Z}}x(n)z^{-n-1},$$
for $x\in {\frak g}$. In the case $k\ne - h^{\vee}$, $N(k,0)$ has a
conformal vector
\begin{equation}
\omega=\frac{1}{2(k+h^{\vee})}\sum_{i=1}^{\dim {\frak g}}
x^{i}(-1)^{2}{\bf 1},
\end{equation}
where $\{x^{i}\}_{i=1, \dots, \dim {\frak g}}$ is an arbitrary
orthonormal basis of ${\frak g}$ with respect to the form $(\cdot,
\cdot)$. We have the following result from \cite{FZ} (see also
\cite{FB}, \cite{K2} \cite{LL}, \cite{L}, \cite{MP}):

\begin{prop}
If $k\ne -h^{\vee}$, the quadruple $(N(k,0), Y, {\bf 1}, \omega)$
defined above is a vertex operator algebra.
\end{prop}

The associative algebra $A(N(k,0))$ is identified in next
proposition:

\begin{prop} [\cite{FZ}, Theorem 3.1.1] \label{p.FZ2}
The associative algebra $A(N(k,0))$
is canonically isomorphic to $U (\frak g ) $.
The isomorphism is given by $F:A(N(k,0)) \to U (\frak g )$
\begin{eqnarray*}
F([x_1(-n_1 -1)\cdots x_m(-n_m -1){\bf 1}])= (-1)^{n_1+\cdots +n_m}
x_m \cdots x_1,
\end{eqnarray*}
for any $x_1, \ldots ,x_m \in \frak g$ and any $n_1, \ldots ,n_m \in {\Z}_{+}$.
\end{prop}

Since every $\hat{\frak g}$-submodule of $N(k,0)$ is also an ideal
in the vertex operator algebra $N(k,0)$, it follows that
$L(k,0)$ is a vertex operator algebra, for any $k\ne -h^{\vee}$.
The associative algebra $A(L(k,0))$ is identified in
the next proposition, in the case when the maximal $\hat{\frak g}$-submodule of
$N(k,0)$ is generated by one singular vector.

\begin{prop} \label{p.1.5.5}
Assume that the maximal $\hat{\frak g}$-submodule of $N(k,0)$ is
generated by a singular vector, i.e. $J(k,0)=U(\hat{\frak g})v.$
Then
$$A(L(k,0)) \cong \frac{U(\frak g)}{I},$$
where $I$ is the two-sided ideal of $U(\frak g)$ generated by
$v'=F([v])$. \\ \noindent Let $U$ be a $\frak g$-module. Then $U$ is
an $A(L(k,0))$-module if and only if $IU=0$.
\end{prop}

\subsection{Modules for associative algebra $A(L(k,0))$}
\label{subsec.3.3}

In this subsection we present the method from \cite{A2}, \cite{AM},
\cite{MP} for classification of irreducible $A(L(k,0))$-modules from
the category $\mathcal{O}$ by solving certain systems of polynomial
equations. We assume that the maximal $\hat{\frak g}$-submodule of
$N(k,0)$ is generated by a singular vector $v$.

Denote by $_L$ the adjoint action of  $U(\frak g)$ on $U(\frak g)$
defined by $ X_Lf=[X,f]$ for $X \in \frak g$ and $f \in U(\frak g)$.
Let $R$ be a $U(\frak g)$-submodule of $U(\frak g)$ generated by the
vector $v'=F([v])$ under the adjoint action. Clearly, $R$ is an
irreducible highest weight $U(\frak g)$-module with the highest
weight vector $v'$. Let $R_{0}$ be the zero-weight subspace of $R$.
The next proposition follows from \cite[Proposition 2.4.1]{A2},
\cite[Lemma 3.4.3]{AM}:

\begin{prop}
Let $V(\mu)$ be an irreducible highest weight $U(\frak g)$-module
with the highest weight vector $v_{\mu}$, for $\mu \in {\frak
h}^{*}$. The following statements are equivalent:
\item[(1)] $V(\mu)$ is an $A(L(k,0))$-module,
\item[(2)] $RV(\mu)=0$,
\item[(3)] $R_{0}v_{\mu}=0.$
\end{prop}

Let $r \in R_{0}$. Clearly there exists the unique polynomial $p_{r}
\in S( \frak h)$ such that
\[
rv_{\mu}=p_{r}(\mu)v_{\mu}.
\]
Set $ {\mathcal P}_{0}=\{ \ p_{r} \ \vert \ r \in R_{0} \}.$ We
have:

\begin{coro} \label{c.1.7.2} There is one-to-one correspondence between
\item[(1)] irreducible $A(L(k,0))$-modules
from the category $\mathcal{O}$,
\item[(2)] weights $\mu \in {\frak h}^{*}$ such that $p(\mu)=0$
for all $p \in {\mathcal P}_{0}$.
\end{coro}

\section{Simple Lie algebra of type $A_{l}$}

Let $\Delta = \{ \epsilon_{i} - \epsilon_{j} \, \vert \ i,j=1,
\ldots ,l+1, \, i \neq j \}$ be the root system of type $A_{l}$. Fix
the set of positive roots $\Delta_{+}= \{ \epsilon_{i} -
\epsilon_{j} \, \vert \ i<j \}$. Then the simple roots are
$\alpha_{1}= \epsilon_{1} - \epsilon_{2}$, $\alpha_{2}= \epsilon_{2}
- \epsilon_{3}$, \ldots , $\alpha_{l}= \epsilon_{l} -
\epsilon_{l+1}$. The highest root is $\theta = \epsilon_{1} -
\epsilon_{l+1}= \alpha_{1}+ \alpha_{2}+ \cdots + \alpha_{l}$.

Let $\frak g$ be the simple Lie algebra associated to the root
system of type $A_{l}$. Let $e_{i},f_{i},h_{i}$, $i=1, \ldots ,l$ be
the Chevalley generators of $\frak g$. Fix the root vectors:
\begin{eqnarray*}
e_{\epsilon_i - \epsilon_j} \!\!\!\!&=& \!\!\![e_{j-1},[e_{j-2},[
\ldots
[e_{i+1},e_{i}] \ldots] \, ] \, ], \qquad i<j, \\
f_{\epsilon_i - \epsilon_j} \!\!\!\!&=& \!\!\! [f_{i},[f_{i+1},[
\ldots [f_{j-2},f_{j-1}] \ldots] \, ] \, ], \qquad i<j.
\end{eqnarray*}
Denote by $h_{\alpha}= \alpha ^{\vee}= [e_{\alpha},f_{\alpha}]$
coroots, for any positive root $\alpha \in \Delta_{+}$. It is clear
that $h_{\alpha_{i}}=h_{i}$. Let ${\frak g}={\frak n}_{-} \oplus
{\frak h} \oplus {\frak n}_{+}$ be the corresponding triangular
decomposition of ${\frak g}$. Denote by $\omega_{1}, \ldots ,
\omega_{l} \in {\frak h}^{*}$ the fundamental weights of $\frak g$,
defined by $\omega_{i}(\alpha ^{\vee} _{j})= \delta _{ij}$ for all
$i,j=1, \ldots ,l$.

\section{Vertex operator algebra
$L(-\frac{1}{2}(l+1),0)$ associated to affine Lie algebra of type
$A_{l}^{(1)}$, for even $l$}

Let $\hat{\frak g}$ be the affine Lie algebra associated to simple
Lie algebra $\frak g$ of type $A_{l}$. We want to show that the
maximal $\hat{\frak g}$-submodule of $N(-\frac{1}{2}(l+1),0)$ is
generated by a singular vector, for even $l$. We need two lemmas to
prove that.

Denote by $\lambda$ the weight $-\frac{1}{2}(l+1) \Lambda_{0}$. Then
$N(-\frac{1}{2}(l+1),0)$ is a quotient of $M(\lambda)$ and
$L(-\frac{1}{2}(l+1),0) \cong L(\lambda)$.

\begin{lem} \label{l.3.2.1}
The weight $\lambda=-\frac{1}{2}(l+1) \Lambda_{0}$ is admissible and
$$ \hat{\Pi}^{\vee}_{\lambda}= \{ (2 \delta - \theta)^{\vee},
\alpha_{1}^{\vee},\alpha_{2}^{\vee}, \ldots , \alpha_{l}^{\vee}
\}.$$
\end{lem}
{\bf Proof:} Clearly
\begin{eqnarray*}
&& \langle \lambda + \rho,\alpha _{i}^{\vee}\rangle = 1 \ \
\mbox{for } i=1,
\ldots ,l, \\
&& \langle \lambda  + \rho,\alpha _{0}^{\vee} \rangle =
-\frac{1}{2}(l-1),
\end{eqnarray*}
which implies
\begin{eqnarray*}
&& \langle \lambda  + \rho,(2 \delta - \theta)^{\vee} \rangle =
\langle \lambda  + \rho, 2\alpha _{0}^{\vee} +\alpha _{1}^{\vee} +
\ldots + \alpha _{l}^{\vee} \rangle = 1.
\end{eqnarray*}
The claim of lemma now follows easily. $\;\;\;\;\Box$

\begin{lem} \label{l.3.2.2}
Vector
$$v= \sum _{i=1}^{l} \frac{l-2i+1}{l+1}h_{i}(-1)e_{\theta}(-1){\bf 1}-
\sum _{i=1}^{l-1} e_{\epsilon_{1} - \epsilon_{i+1}}(-1)
e_{\epsilon_{i+1} - \epsilon_{l+1}}(-1){\bf 1} - \frac{1}{2}(l-1)
e_{\theta}(-2){\bf 1} $$ is a singular vector in
$N(-\frac{1}{2}(l+1),0)$.
\end{lem}
{\bf Proof:} It can be directly verified that
\begin{eqnarray*}
& & e_{j}(0).v=0, \ j=1, \ldots ,l, \\
& & f_{\theta}(1).v=0. \  \;\;\;\;\Box
\end{eqnarray*}

\begin{thm} \label{t.3.3} The maximal $\hat{\frak g}$-submodule
of $N(-\frac{1}{2}(l+1),0)$ is $J(-\frac{1}{2}(l+ \nolinebreak
1),0)= \linebreak U(\hat{\frak g})v$, where
$$v= \sum _{i=1}^{l} \frac{l-2i+1}{l+1}h_{i}(-1)e_{\theta}(-1){\bf 1}- \sum _{i=1}^{l-1}
e_{\epsilon_{1} - \epsilon_{i+1}}(-1) e_{\epsilon_{i+1} -
\epsilon_{l+1}}(-1){\bf 1} - \frac{1}{2}(l-1) e_{\theta}(-2){\bf 1}.
$$
\end{thm}
{\bf Proof:} It follows from Theorem \ref{t.KW1} and Lemma
\ref{l.3.2.1} that the maximal submodule of the Verma module
$M(\lambda)$ is generated by $l+1$ singular vectors with weights
$$r_{2 \delta - \theta}.\lambda, r_{\alpha_{1}}.\lambda,
\ldots ,r_{\alpha_{l}}.\lambda.$$ It follows from Lemma
\ref{l.3.2.2} that $v$ is a singular vector of weight $\lambda-2
\delta + \theta=r_{2 \delta - \theta}.\lambda$. Other singular
vectors have weights
$$r_{\alpha_{i}}.\lambda=\lambda- \langle \lambda+
\rho, \alpha_{i}^{\vee} \rangle \alpha_{i}=\lambda - \alpha_{i}, \
i=1,\ldots ,l,$$ so the images of these vectors under the projection
of $M(\lambda)$ onto $N(-\frac{1}{2}(l+1),0)$ are 0. Therefore, the
maximal submodule of $N(-\frac{1}{2}(l+1),0)$ is generated by the
vector $v$, i.e. $J(-\frac{1}{2}(l+1),0)= U(\hat{\frak g})v$.
$\;\;\;\;\Box$

It follows that
$$L(-\frac{1}{2}(l+1),0) \cong \frac{N(-\frac{1}{2}(l+1),0)}{U(\hat{\frak
g})v}.$$ Using Theorem \ref{t.3.3} and Proposition \ref{p.1.5.5} we
can determine the associative algebra $A(L(-\frac{1}{2}(l+1),0))$:

\begin{prop} \label{t.3.2.5}
The associative algebra $A(L(-\frac{1}{2}(l+1),0))$ is isomorphic to
the algebra $U(\frak g)/I$, where $I$ is the two-sided ideal of
$U(\frak g)$ generated by
$$v'= \sum _{i=1}^{l} \frac{l-2i+1}{l+1}h_{i}e_{\theta}-
\sum _{i=1}^{l-1} e_{\epsilon_{i+1} - \epsilon_{l+1}}
e_{\epsilon_{1} - \epsilon_{i+1}} + \frac{1}{2}(l-1) e_{\theta}. $$
\end{prop}
{\bf Proof:} The maximal submodule $J(-\frac{1}{2}(l+1),0)$ od
$N(-\frac{1}{2}(l+1),0)$ is generated by the singular vector $v$. It
follows from Proposition \ref{p.1.5.5} that
$$A(L(-\frac{1}{2}(l+1),0)) \cong \frac{U(\frak g)}{I},$$
where $I$ is the two-sided ideal in $U(\frak g)$ generated by
$v'=F([v])$. Proposition \ref{p.FZ2} now implies that
\begin{eqnarray*}
&&v'= F([v])=\sum _{i=1}^{l} \frac{l-2i+1}{l+1}e_{\theta}h_{i}- \sum
_{i=1}^{l-1} e_{\epsilon_{i+1} - \epsilon_{l+1}}
 e_{\epsilon_{1} - \epsilon_{i+1}} + \frac{1}{2}(l-1)
e_{\theta} \\
&& \quad = \sum _{i=1}^{l} \frac{l-2i+1}{l+1}h_{i}e_{\theta}- \sum
_{i=1}^{l-1} e_{\epsilon_{i+1} - \epsilon_{l+1}}
 e_{\epsilon_{1} - \epsilon_{i+1}} + \frac{1}{2}(l-1)
e_{\theta},
\end{eqnarray*}
which implies the claim of proposition. $\;\;\;\;\Box$

\section{Classification of irreducible weak
$L(-\frac{1}{2}(l+1),0)$-modules from category $\mathcal{O}$}

In this section we classify irreducible weak
$L(-\frac{1}{2}(l+1),0)$-modules that are in category $\mathcal{O}$
as $\hat{\frak g}$-modules, using methods from  \cite{A2},
\cite{AM}, \cite{MP} presented in Subsection \ref{subsec.3.3}.
First, we determine a basis for the vector space ${\mathcal P}_{0}$
defined in that subsection. Recall that $_L$ denotes the adjoint
action of $U(\frak g)$ on $U(\frak g)$ defined by $ X_Lf=[X,f]$ for
$X \in \frak g$ and $f \in U(\frak g)$.

\begin{lem} \label{l.3.3.1}  Let
\begin{eqnarray*}
&\!\!\!\!\!\! &\!\! p_{i}(h)=h_{i} \left( \sum _{j=1}^{i-1}
\frac{-2j}{l+1}h_{j}+ \frac{l-2i+1}{l+1}h_{i}+ \sum _{j=i+1}^{l}
\frac{2l-2j+2}{l+1}h_{j} + \frac{1}{2}(l+1)-i \right),
\end{eqnarray*}
for $i=1, \ldots ,l$. Then $p_{1}, \ldots ,p_{l} \in {\mathcal
P}_{0}$.
\end{lem}
{\bf Proof:} We claim that
\begin{eqnarray} \label{rel.polinomi}
&&(-1)^{i}(f_{i}f_{i-1} \ldots f_{1}f_{i+1} \ldots f_{l})_L \, v'
\in p_{i}(h)  + U(\frak g){\frak n}_{+}, \ \mbox{for } i=1, \ldots
,l. \nonumber \\
&& \mbox{}
\end{eqnarray}
One can easily verify that for $i \in \{ 1, \ldots ,n \}$ the
following relations hold:
\begin{eqnarray*}
&&(f_{i}f_{i-1} \ldots f_{1}f_{i+1} \ldots f_{l})_L \, e_{\theta}
=(-1)^{i}h_{i}, \\
&&(f_{i}f_{i-1} \ldots f_{1})_L \, e_{\epsilon_{1} - \epsilon_{j+1}}
=(-1)^{i}f_{\epsilon_{j+1} - \epsilon_{i+1}}, \quad j<i, \\
&&(f_{i+1} \ldots f_{l})_L \, e_{\epsilon_{j+1} - \epsilon_{l+1}}
=e_{\epsilon_{j+1} - \epsilon_{i+1}}, \quad j<i, \\
&&(f_{i-1}f_{i-2} \ldots f_{1})_L \, e_{\epsilon_{1} -
\epsilon_{j+1}}
=(-1)^{i-1}f_{\epsilon_{j+1} - \epsilon_{i}}, \quad j<i-1, \\
&&(f_{i}f_{i+1} \ldots f_{l})_L \, e_{\epsilon_{j+1} -
\epsilon_{l+1}} =e_{\epsilon_{j+1} - \epsilon_{i}}, \quad j<i-1.
\end{eqnarray*}
Using Proposition \ref{t.3.2.5}, we obtain
\begin{eqnarray*} && (f_{i}f_{i-1} \ldots f_{1}f_{i+1}
\ldots f_{l})_L \, v' \in
(-1)^{i} \sum _{j=1}^{l} \frac{l-2j+1}{l+1}h_{j}h_{i} \\
&& \quad - (-1)^{i} \sum _{j=1}^{i-2} ( h_{i}h_{j}+
e_{\epsilon_{j+1} - \epsilon_{i+1}} f_{\epsilon_{j+1} -
\epsilon_{i+1}} - e_{\epsilon_{j+1} - \epsilon_{i}}
f_{\epsilon_{j+1} - \epsilon_{i}}) - (-1)^{i}( h_{i}h_{i-1} +
e_{i}f_{i})  \\
&& \quad -(-1)^{i+1} \sum _{j=i+1}^{l}h_{j}h_{i} +
\frac{1}{2}(l-1)(-1)^{i}h_{i} +  U(\frak g){\frak n}_{+}.
\end{eqnarray*}
Since
\begin{eqnarray*}
&& e_{\epsilon_{j+1} - \epsilon_{i+1}} f_{\epsilon_{j+1} -
\epsilon_{i+1}} - e_{\epsilon_{j+1} - \epsilon_{i}}
f_{\epsilon_{j+1} - \epsilon_{i}} \in h_{\epsilon_{j+1} -
\epsilon_{i+1}}- h_{\epsilon_{j+1} - \epsilon_{i}} +  U(\frak
g){\frak n}_{+} = h_{i} +  U(\frak g){\frak n}_{+},
\end{eqnarray*}
it follows that
\begin{eqnarray*}
&& (-1)^{i} (f_{i}f_{i-1} \ldots f_{1}f_{i+1} \ldots f_{l})_L \, v'
\in
\sum _{j=1}^{l} \frac{l-2j+1}{l+1}h_{j}h_{i} \\
&& \quad - \sum _{j=1}^{i-1}h_{j}h_{i} + \sum _{j=i+1}^{l}h_{j}h_{i}
-(i-1)h_{i} + \frac{1}{2}(l-1)h_{i} +  U(\frak g){\frak n}_{+},
\end{eqnarray*}
which implies relation (\ref{rel.polinomi}). $\;\;\;\;\Box$

\begin{lem} \label{l.3.3.2}
\begin{eqnarray*}
{\mathcal P}_{0}= \mbox{span}_{\C} \{ p_{1}, \ldots ,p_{l} \}.
\end{eqnarray*}
\end{lem}
{\bf Proof:} Lemma \ref{l.3.3.1} implies that $p_{1}, \ldots ,p_{l}$
are linearly independent polynomials in the set ${\mathcal P}_{0}$.
It follows from the definition of set ${\mathcal P}_{0}$ that $\dim
{\mathcal P}_{0} = \dim R_{0}$, where $R$ is the highest weight
$U(\frak g)$-module with highest weight $\theta$, and $R_{0}$ the
zero-weight subspace of $R$. Since $R$ is isomorphic to the adjoint
module for $\frak g$, it follows that $\dim R_{0}=l$. Thus,
polynomials $p_{1}, \ldots ,p_{l}$ form a basis for ${\mathcal
P}_{0}$. $\;\;\;\;\Box$

\begin{prop} \label{t.3.3.3}
For every subset $S=\{ i_{1}, \ldots ,i_{k} \} \subseteq \{1,2,
\ldots , l \}$, $i_{1}< \ldots <i_{k}$, we define weights
\begin{eqnarray*}
&& \mu _{S}= \sum _{j=1}^{k}\left( \sum _{s=j+1}^{k} (-1)^{s-j}i_{s}
+ \sum _{s=1}^{j-1}(-1)^{j-s+1}i_{s}+ (-1)^{k-j+1}\frac{l+1}{2}
\right) \omega _{i_{j}},
\end{eqnarray*}
where $\omega _{1}, \ldots , \omega _{l}$ are fundamental weights
for $\frak g$. Then the set
$$ \{ V(\mu _{S}) \ \vert \ S \subseteq \{1,2, \ldots , l \} \}$$
provides the complete list of irreducible
$A(L(-\frac{1}{2}(l+1),0))$-modules from the category $\mathcal{O}$.
\end{prop}
{\bf Proof:} It follows from Corollary \ref{c.1.7.2} and Lemma
\ref{l.3.3.2} that highest weights $\mu \in {\frak h}^{*}$ of
irreducible $A(L(-\frac{1}{2}(l+1),0))$-modules $V(\mu)$ are in
one-to-one correspondence with solutions of the system of polynomial
equations
\begin{eqnarray}
&&\!\!\!\!\!\!\! h_{i} \left( -\sum _{j=1}^{i-1} jh_{j}+
\frac{l-2i+1}{2}h_{i}+ \sum _{j=i+1}^{l} (l-j+1)h_{j}
+ \frac{1}{4}(l+1)^{2}- \frac{1}{2}i(l+1) \right)=0, \nonumber \\
&& \mbox{} \label{3.3.1.1}
\end{eqnarray}
for $i=1, \ldots ,l$.

Let $i,j \in \{1,2, \ldots , l \}$, $i<j$. If we multiply the $i$-th
equation of system (\ref{3.3.1.1}) by $h_{j}$, and the $j$-th
equation by $h_{i}$ and then subtract these equations, we obtain
\begin{eqnarray} \label{3.3.1}
&& h_{i}h_{j} \left(h_{i}+2h_{i+1}+2h_{i+2}+ \ldots +
2h_{j-1}+h_{j}+j-i \right)=0.
\end{eqnarray}

Let $S=\{ i_{1}, \ldots ,i_{k} \}$, $i_{1}< \ldots <i_{k}$ be the
subset of $\{1,2, \ldots , l \}$ such that $h_{i}=0$ for $i \notin
S$ and $h_{i} \neq 0$ for $i \in S$. From relation (\ref{3.3.1}) and
relation (\ref{3.3.1.1}) for $i=i_{k}$, we get the system
\begin{eqnarray}
& h_{i_{1}}+h_{i_{2}}+i_{2}-i_{1} =0& \nno \\
& h_{i_{2}}+h_{i_{3}}+i_{3}-i_{2} =0& \nno \\
& \qquad \qquad \vdots &   \label{3.3.3} \\
& h_{i_{k-1}}+h_{i_{k}}+i_{k}-i_{k-1} =0& \nno \\
& -i_{1}h_{i_{1}}-i_{2}h_{i_{2}}- \ldots - i_{k-1}h_{i_{k-1}}+
\frac{l-2i_{k}+1}{2}h_{i_{k}} + \frac{1}{4}(l+1)^{2}-
\frac{1}{2}i_{k}(l+1)=0.& \nno
\end{eqnarray}
If we multiply the first equation of system (\ref{3.3.3}) by
$i_{1}$, the second equation by $i_{2}-i_{1}$, the third equation by
$i_{3}-i_{2}+i_{1}$, \ldots , the $(k-1)$-th equation by
$i_{k-1}-i_{k-2}+ \ldots +(-1)^{k}i_{1}$ and then sum these
equations and the $k$-th equation, we obtain:
\begin{eqnarray*}
 \left( \frac{l+1}{2}- i_{k}+i_{k-1}- \ldots +(-1)^{k}i_{1}
\right) \left(h_{i_{k}}+\frac{l+1}{2}- i_{k-1}+i_{k-2}- \ldots
+(-1)^{k-1}i_{1} \right)=0.
\end{eqnarray*}
Since $l$ is even, we have
\begin{eqnarray*}
 \frac{l+1}{2}- i_{k}+i_{k-1}- \ldots +(-1)^{k}i_{1} \neq 0,
\end{eqnarray*}
which implies
\begin{eqnarray*}
h_{i_{k}}=i_{k-1}-i_{k-2}+ \ldots +(-1)^{k}i_{1}-\frac{l+1}{2}.
\end{eqnarray*}
Using the first $k-1$ equations of system (\ref{3.3.3}) one can
easily obtain that
\begin{eqnarray*}
&& h_{i_{j}}=  \sum _{s=j+1}^{k} (-1)^{s-j}i_{s} + \sum
_{s=1}^{j-1}(-1)^{j-s+1}i_{s}+ (-1)^{k-j+1}\frac{l+1}{2}, \  j=1,
\ldots ,k
\end{eqnarray*}
is a solution of this system. Thus, $V(\mu _{S})$ is an irreducible
$A(L(-\frac{1}{2}(l+1),0))$-module which implies the claim of
proposition. $\;\;\;\;\Box$

It follows from Zhu's theory that:

\begin{thm} \label{c.3.3.4}
The set
$$ \{ L(-\frac{1}{2}(l+1), \mu _{S}) \ \vert \
S \subseteq \{1,2, \ldots , l \} \}$$ provides the complete list of
irreducible weak $L(-\frac{1}{2}(l+1),0)$-modules from the category
$\mathcal{O}$.
\end{thm}

Theorem \ref{c.3.3.4} implies that there are $2^{l}$ irreducible
weak $L(-\frac{1}{2}(l+1),0)$-modules from category $\mathcal{O}$.
The weight $\mu _{S}$ is a dominant integral weight for $\frak g$ if
and only if $S=\emptyset$, i.e. if and only if $\mu _{S}=0$. It
follows that

\begin{coro} \label{c.3.3.4.1}
$L(-\frac{1}{2}(l+1),0)$ is the only irreducible
$L(-\frac{1}{2}(l+1),0)$-module.
\end{coro}

\section{Complete reducibility of weak
$L(-\frac{1}{2}(l+1),0)$-modules from category $\mathcal{O}$}
\label{sec.3.4}

In this section we show that every weak
$L(-\frac{1}{2}(l+1),0)$-module from category $\mathcal{O}$ is
completely reducible. We introduce the notation
$\lambda_{S}=-\frac{1}{2}(l+1) \Lambda_{0}+ \mu _{S}$, for every $S
\subseteq \{1,2, \ldots , l \}$. The following lemma is crucial for
proving complete reducibility.

\begin{lem}\label{l.3.3.5}
The weight $\lambda _{S} \in \hat{\frak h} ^{*}$ is admissible, for
every $S \subseteq \{1,2, \ldots , l \}$.
\end{lem}
{\bf Proof:} We have to show
\begin{eqnarray}
& &\langle \lambda_{S} + \rho,\tilde{\alpha}^{\vee}\rangle \notin
-{\Z}_{+} \mbox{ for all }
\tilde{\alpha} \in \hat{\Delta}^{\mbox{\scriptsize{re}}}_{+}, \label{3.3.4}\\
& &{\Q} \hat{\Delta}^{\vee \mbox{\scriptsize{re}}}_{\lambda
_{S}}={\Q} \hat{\Pi}^{\vee} .\label{3.3.5}
\end{eqnarray}
First, let us prove relation (\ref{3.3.4}). Any positive real root
$\tilde{\alpha} \in \hat{\Delta}^{\mbox{\scriptsize{re}}}_{+}$ of
$\hat{\frak g}$ is of the form $\tilde{\alpha}= \alpha +m \delta$,
for $m>0$ and $\alpha \in \Delta$ or $m=0$ and $\alpha \in
{\Delta}_{+}$. Positive roots of $\frak g$ are $\epsilon_{i} -
\epsilon_{j}$, $i<j$, and negative roots are $-(\epsilon_{i} -
\epsilon_{j})$, $i<j$.

Clearly, $( \bar{\rho}, \epsilon_{i} - \epsilon_{j})=j-i$. Let $s,t
\in \{ 1, \ldots ,k \}$ be the indices such that $S \cap \{ i,i+1,
\ldots ,j-1 \}= \{ i_{s}, \ldots ,i_{t} \}$. Clearly, $i_{s} \geq i$
and $i_{t} \leq j-1$. Furthermore, $( \mu _{S}, \epsilon_{i} -
\epsilon_{j})=h_{i_{s}}+ \ldots \ +h_{i_{t}}$.

We obtain
\begin{eqnarray}
& & \langle \lambda_{S} + \rho,\tilde{\alpha}^{\vee}\rangle =
\frac{1}{2}m(l+1)+( \bar{\rho} , \alpha )+ ( \mu _{S} , \alpha ) ,
\label{3.3.6}
\end{eqnarray}
where $\bar{\rho}$ is the sum of fundamental weights of $\frak g$.
Let $S=\{ i_{1}, \ldots ,i_{k} \} \subseteq \{1,2, \ldots , l \}$,
$i_{1}< \ldots <i_{k}$. Proposition \ref{t.3.3.3} implies that
$\mu_{S}= \sum _{j=1}^{k} h_{i_{j}} \omega _{i_{j}}$, where
\begin{eqnarray*} && h_{i_{j}}=  \sum _{s=j+1}^{k}
(-1)^{s-j}i_{s} + \sum _{s=1}^{j-1}(-1)^{j-s+1}i_{s}+
(-1)^{k-j+1}\frac{l+1}{2}, \ j=1, \ldots ,k.
\end{eqnarray*}

First consider the case $\alpha= \epsilon_{i} - \epsilon_{j}$, $i<j$
and $m \geq 0$.

If $t-s+1$ is even, then using relations from system (\ref{3.3.3})
we get
\begin{eqnarray*}
( \mu _{S}, \epsilon_{i} - \epsilon_{j})= -(i_{s+1}-i_{s})- \ldots -
(i_{t}-i_{t-1}) \geq -(i_{t}-i_{s}),
\end{eqnarray*}
and relation (\ref{3.3.6}) implies
\begin{eqnarray*}
&& \langle \lambda_{S} + \rho,\tilde{\alpha}^{\vee}\rangle \geq (
\bar{\rho} , \epsilon_{i} - \epsilon_{j} )+ ( \mu _{S} ,
\epsilon_{i} - \epsilon_{j} )
\geq (j-i)- (i_{t}-i_{s}) \\
&& \qquad \qquad \quad \ \ = (j-i_{t})+ (i_{s}-i) >0.
\end{eqnarray*}

Suppose now that $t-s+1$ is odd. Then $( \mu _{S}, \epsilon_{i} -
\epsilon_{j}) \notin \Z$, and if $m=0$, then $\langle \lambda_{S} +
\rho,\tilde{\alpha}^{\vee}\rangle \notin \Z$. Let $m \geq 1$. Then
\begin{eqnarray*}
&& ( \mu _{S}, \epsilon_{i} - \epsilon_{j})=h_{i_{s}}+ \ldots \
+h_{i_{t-1}} +h_{i_{t}}= -(i_{s+1}-i_{s})- \ldots -
(i_{t-1}-i_{t-2})+h_{i_{t}} \\
&& \qquad \qquad \quad \  \geq -(i_{t-1}-i_{s})+h_{i_{t}}.
\end{eqnarray*}
We have
\begin{eqnarray*}
&& h_{i_{t}}=  \sum _{s=t+1}^{k} (-1)^{s-t}i_{s} + \sum
_{s=1}^{t-1}(-1)^{t-s+1}i_{s}+ (-1)^{k-t+1}\frac{l+1}{2}.
\end{eqnarray*}
Clearly, $\sum _{s=1}^{t-1}(-1)^{t-s+1}i_{s} \geq 0$. If $k-t$ is
even, then
\begin{eqnarray*}
h_{i_{t}} \geq (i_{t+2}-i_{t+1})+ \ldots + (i_{k}-i_{k-1}) -
\frac{1}{2}(l+1) \geq - \frac{1}{2}(l+1),
\end{eqnarray*}
which implies
\begin{eqnarray*}
( \mu _{S}, \epsilon_{i} - \epsilon_{j}) \geq -(i_{t-1}-i_{s})-
\frac{1}{2}(l+1).
\end{eqnarray*}
It follows that
\begin{eqnarray*}
&& \langle \lambda_{S} + \rho,\tilde{\alpha}^{\vee}\rangle \geq
\frac{1}{2}(l+1) +( \bar{\rho} , \epsilon_{i} - \epsilon_{j} )+ (
\mu _{S} ,
\epsilon_{i} - \epsilon_{j}) \\
&& \qquad \qquad \quad \ \geq \frac{1}{2}(l+1)+ (j-i)-(i_{t-1}-i_{s})-\frac{1}{2}(l+1) \\
&& \qquad \qquad \quad \ = (j-i_{t-1})+ (i_{s}-i) >0.
\end{eqnarray*}
If $k-t$ is odd, then
\begin{eqnarray*}
h_{i_{t}} \geq (i_{t+2}-i_{t+1})+ \ldots + (i_{k-1}-i_{k-2})- i_{k}
+\frac{1}{2}(l+1) \geq \frac{1}{2}(l+1) - i_{k},
\end{eqnarray*}
which implies
\begin{eqnarray*}
( \mu _{S}, \epsilon_{i} - \epsilon_{j}) \geq -(i_{t-1}-i_{s})+
\frac{1}{2}(l+1)- i_{k}.
\end{eqnarray*}
We obtain
\begin{eqnarray*}
&& \langle \lambda_{S} + \rho,\tilde{\alpha}^{\vee}\rangle \geq
\frac{1}{2}(l+1) +( \bar{\rho} , \epsilon_{i} - \epsilon_{j} )+ (
\mu _{S} , \epsilon_{i} - \epsilon_{j} )
\geq \frac{1}{2}(l+1)+ (j-i) \\
&& \qquad  \quad \ -(i_{t-1}-i_{s})+ \frac{1}{2}(l+1) -i_{k} =
(l-i_{k})+(j-i_{t-1})+ (i_{s}-i)+1 >0.
\end{eqnarray*}

Thus, we have proved that, if $\alpha= \epsilon_{i} - \epsilon_{j}$,
$i<j$ and $m \geq 0$, then $ \langle \lambda_{S} +
\rho,\tilde{\alpha}^{\vee}\rangle \notin -{\Z}_{+}$.

Now, let us consider the case $\alpha= -(\epsilon_{i} -
\epsilon_{j})$, $i<j$ and $m \geq 1$.

Then
\begin{eqnarray*}
& & \langle \lambda_{S} + \rho,\tilde{\alpha}^{\vee}\rangle =
\frac{1}{2}m(l+1)-( \bar{\rho}, \epsilon_{i} - \epsilon_{j} )- ( \mu
_{S} , \epsilon_{i} - \epsilon_{j} ).
\end{eqnarray*}

If $t-s+1$ is even, then $( \mu _{S}, \epsilon_{i} - \epsilon_{j})$
is an integer and $( \mu _{S}, \epsilon_{i} - \epsilon_{j}) \leq 0$,
so if $m$ is odd, then $\langle \lambda_{S} +
\rho,\tilde{\alpha}^{\vee}\rangle \notin \Z$. If $m$ is even, then
$m \geq 2$, and we get
\begin{eqnarray*}
\langle \lambda_{S} + \rho,\tilde{\alpha}^{\vee}\rangle \geq
(l+1)-(j-i)=(l-j)+i+1>0.
\end{eqnarray*}
If $t-s+1$ is odd, then
\begin{eqnarray*}
( \mu _{S}, \epsilon_{i} - \epsilon_{j})=h_{i_{s}}+ \ldots \
+h_{i_{t-1}} +h_{i_{t}} =-(i_{s+1}-i_{s})- \ldots -
(i_{t-1}-i_{t-2})+ h_{i_{t}},
\end{eqnarray*}
which implies
\begin{eqnarray*}
&& \langle \lambda_{S} + \rho,\tilde{\alpha}^{\vee}\rangle \geq
\frac{1}{2}(l+1)-( j-i )+(i_{s+1}-i_{s})+ \ldots +
(i_{t-1}-i_{t-2})- h_{i_{t}}.
\end{eqnarray*}
If $k-t$ is even, then
\begin{eqnarray*}
h_{i_{t}} = i_{t-1}-i_{t-2}+ \ldots + (-1)^{t-1} i_{2} +
(-1)^{t}i_{1}+i_{k}- i_{k-1}+  \ldots + i_{t+2}- i_{t+1}-
\frac{1}{2}(l+1),
\end{eqnarray*}
which implies
\begin{eqnarray*}
&& \langle \lambda_{S} + \rho,\tilde{\alpha}^{\vee}\rangle \geq
\frac{1}{2}(l+1)-( j-i )+(i_{s-1}-i_{s-2}+ \ldots +
(-1)^{t+1}i_{1}) \\
&& \qquad \qquad  \quad -i_{k}+ (i_{k-1}- i_{k-2})+ \ldots +
(i_{t+3}- i_{t+2})+i_{t+1}+ \frac{1}{2}(l+1).
\end{eqnarray*}
Clearly
\begin{eqnarray*}
i_{s-1}-i_{s-2}+ \ldots + (-1)^{t+1}i_{1} \geq 0,
\end{eqnarray*}
from which we get
\begin{eqnarray*}
\langle \lambda_{S} + \rho,\tilde{\alpha}^{\vee}\rangle \geq
(l-i_{k})+i + (i_{t+1}-j)+ (i_{k-1}- i_{k-2})+ \ldots + (i_{t+3}-
i_{t+2})+1 >0.
\end{eqnarray*}
If $k-t$ is odd, then
\begin{eqnarray*}
h_{i_{t}} = i_{t-1}-i_{t-2}+ \ldots + (-1)^{t}i_{1}-i_{k}+ i_{k-1}-
\ldots + i_{t+2}- i_{t+1}+ \frac{1}{2}(l+1),
\end{eqnarray*}
which implies
\begin{eqnarray*}
&& \langle \lambda_{S} + \rho,\tilde{\alpha}^{\vee}\rangle \geq
\frac{1}{2}(l+1)-( j-i )+(i_{s-1}-i_{s-2}+ \ldots +
(-1)^{t+1}i_{1}) \\
&& \qquad \qquad  \quad
+i_{k}-i_{k-1}+ \ldots + i_{t+3}- i_{t+2}+i_{t+1}- \frac{1}{2}(l+1)  \\
&& \qquad \qquad  \quad \geq (i_{t+1}-j)+i+ (i_{k}-i_{k-1})+ \ldots
+ (i_{t+3}- i_{t+2})>0.
\end{eqnarray*}
We have proved that, if $\alpha=- (\epsilon_{i} - \epsilon_{j})$,
$i<j$ and $m \geq 1$, then $ \langle \lambda_{S} +
\rho,\tilde{\alpha}^{\vee}\rangle \notin -{\Z}_{+}$.

Thus, we have verified the relation (\ref{3.3.4}). Moreover, one can
easily check that coroots
\begin{eqnarray*}
&& (\delta - \alpha _{i_{j}})^{\vee}, \ j=1, \ldots , k, \\
&& \alpha _{i_{j}}^{\vee}+\alpha _{i_{j}+1}^{\vee}+ \ldots + \alpha
_{i_{j+1}}^{\vee},
\ j=1, \ldots , k-1, \\
&& \alpha _{i}^{\vee}, \ i \notin S, \ i \in \{1,2, \ldots ,l \}
\end{eqnarray*}
are elements of the set $\hat{\Delta}^{\vee
\mbox{\scriptsize{re}}}_{\lambda _{S}}$ which implies ${\Q}
\hat{\Delta}^{\vee \mbox{\scriptsize{re}}}_{\lambda _{S}}={\Q}
\hat{\Pi}^{\vee}$, and relation (\ref{3.3.5}) is also proved.
$\;\;\;\;\Box$

\begin{thm} \label{t.3.3.6}
Let $M$ be a weak $L(-\frac{1}{2}(l+1),0)$-module from the category
$\mathcal{O}$. Then $M$ is completely reducible.
\end{thm}
{\bf Proof:} Let $L(\lambda)$ be some irreducible subquotient of
$M$. Then $L(\lambda)$ is an irreducible weak
$L(-\frac{1}{2}(l+1),0)$-module, and Theorem \ref{c.3.3.4} implies
that there exists $S \subseteq \{1,2, \ldots , l \}$ such that
$\lambda=-\frac{1}{2}(l+1) \Lambda_{0}+ \mu _{S}$. It follows from
Lemma \ref{l.3.3.5} that such $\lambda$ is admissible. Theorem
\ref{t.KW2} now implies that $M$ is completely reducible.
$\;\;\;\;\Box$

\begin{rem} Using the fact that $L(-\frac{1}{2}(l+1),0)$ is the only
irreducible $L(-\frac{1}{2}(l+1),0)$-module (Corollary
\ref{c.3.3.4.1}), one can show (as in \cite[Lemma 26]{P1}) that
every $L(-\frac{1}{2}(l+1),0)$-module is in the category
$\mathcal{O}$ as $\hat{\frak g}$-module. It follows now from Theorem
\ref{t.3.3.6} that every $L(-\frac{1}{2}(l+1),0)$-module is
completely reducible, which implies that it is a direct sum of
copies of $L(-\frac{1}{2}(l+1),0)$.
\end{rem}

\bibliography{thesis}
\bibliographystyle{plain}

\vskip 1cm

Department of Mathematics, University of Zagreb, Bijeni\v{c}ka 30,
\linebreak 10000 Zagreb, Croatia

E-mail address: perse@math.hr

\end{document}